\documentclass[12pt]{amsart}
\usepackage{amsmath,amssymb,amsthm}

\usepackage{setspace}
\usepackage[colorlinks=true,linkcolor=blue,citecolor=blue]{hyperref}

\usepackage[cmtip,all]{xy}

\usepackage{xcolor}
\newtheorem{theorem}{Theorem}[section]

\newtheorem{defn}{Definition}[section]
\newtheorem{lemma}[theorem]{Lemma}
\newtheorem{prop}[theorem]{Proposition}

\newtheorem{cor}[theorem]{Corollary}

\theoremstyle{remark}
\newtheorem{rem}[theorem]{Remark}

\newtheorem*{thma}{\bf{Theorem A}}

\newtheorem*{conja}{{\bf Conjecture}}

\newcommand{\cP}{\mathcal{P}}
\newcommand{\cS}{\mathcal{S}}

\newcommand{\cC}{\mathcal{C}}
\newcommand{\cT}{\mathcal{T}}

\begin{document}

\title
{Sums and products in Dynamically Central Syndetic Sets}

\author{Arpita Ghosh}
\curraddr{Department of Mathematics, University of Haifa, Israel}
\email{arpi.arpi16@gmail.com}

\subjclass[2020]{Primary: 05D10. Secondary: 37B20.}

\keywords{Hindman's Conjecture, Sums and products, Dynamically central syndetic sets}

\begin{abstract}
    Glasscock and Le introduced the notion \textit{dynamically central syndetic sets} in a minimal topological dynamical system under the minimal action of the semigroup $(\mathbb{N}, +).$  
    In this short note, we introduce these large sets in the semigroup $(\mathbb{N}, \cdot)$ and establish that dynamically central syndetic sets contain the Hindman configuration $\{x,y,x+y,xy\}$. Furthermore, we show that there are dynamically central piecewise syndetic sets that do not contain such configurations.
\end{abstract}
\maketitle

\section{Introduction}

The interaction between addition and multiplication is of particular interest in Ramsey theory. While classical results in Ramsey theory guarantee the existence of monochromatic additive or multiplicative configurations under finite colourings, the simultaneous occurrence of both types of operation often requires substantially stronger arguments. In this connection, the \textit{Hindman's Conjecture} is one of cornerstones, and several variants of it have been the subject of recent attention in \cite{ABS25, B25, BS24, BM17, BM18, GG, GS16, M17, P21, S18}.

 \begin{conja}[Hindman \cite{H79}] Any finite coloring of $\mathbb{N}$ contains monochromatic sets of the form $\{x, y, xy, x+y\}.$
 \end{conja}

 This problem illustrates the difficulty of understanding the relationship between additive and multiplicative structures in the natural numbers.

A different approach to studying the existence of such configurations is through notions of combinatorial largeness arising from topological dynamics. In this direction, Glasscock and Le introduced and developed the notion of dynamically central syndetic sets in their work \cite{GL25, GL26} and the Combinatorics of \emph{Syndetic, Idempotent Filters.} Their work establishes a connection between topological dynamics and the algebraic structure of families of subsets of $(\mathbb{N}, +)$. In particular, they characterize dynamically central syndetic sets in terms of membership in syndetic, idempotent filters. This characterization provides a useful framework for understanding the combinatorial properties of these sets. In this regard, recall a few definitions as follows.

\begin{defn}
\begin{enumerate}
    \item Let $(X,T)$ be a \emph{minimal topological dynamical system} if $X$ is a compact metric space and $T\colon X\to X$ is a continuous map such that every forward orbit is dense, that is , $\overline{\{T^nx: n \ge 1\}}\cong X$ for all $x \in X$.

   \item  For $x\in X$ and a nonempty open set $U\subseteq X$, the \emph{return-time set} is
\[
R(x,U):=\{n\in\mathbb N:\ T^nx\in U\}.
\]

\item  A set $A\subseteq\mathbb N$ is \emph{dynamically syndetic} if $A\supseteq R(x,U)$ for some minimal system $(X,T)$, some $x\in X$ and some nonempty open $U$. 

\item If moreover $x\in U$, then $A$ is \emph{dynamically central syndetic}. We write $d\mathcal S$ and $dc\mathcal S$ for these two families.
\end{enumerate}
\end{defn}

 Thus, dynamically central syndetic sets inherit a strong form of combinatorial richness from the underlying dynamical structure. This makes them natural objects for investigating the existence of arithmetic configurations. We denote the \textit{family} of dynamically (central) syndetic sets as $d\mathcal S$ $(dc\mathcal{S}).$

 The \textit{dual} of a family $\mathcal{F}$ of subsets of $\mathbb{N}$ is 
 \[
 \mathcal{F}^* := \{B \subseteq \mathbb{N} \mid  \forall A \in \mathcal{F}, B \cap A \neq \emptyset\}
 \]

Dual of the family of dynamically syndetic sets are the family  of \textit{dynamically thick} $(d\mathcal{T})$ and dual of the family of dynamically central syndetic sets are the family  of \textit{dynamically central thick} $(dc\mathcal{T})$ sets. \textit{Thick sets} are the subsets of $\mathbb{N}$ that contain arbitrary long intervals, are the most immediate examples of dynamically thick sets. The family of sets formed by intersections of dynamically (central) syndetic and dynamically (central) thick sets are called \textit{dynamically piecewise syndetic} (\textit{dynamically central piecewise syndetic}) sets. We write $d\mathcal{PS}$ and $dc\mathcal{PS}$ for these two families, respectively. For further examples of these sets are given in \cite{GL25, GL26}.

Motivated by these developments, in multiplicative setting,  a \emph{minimal multiplicative system} is a compact metric space $X$ with a minimal action $n\mapsto T_n$ of $(\mathbb N,\cdot)$ by continuous maps $T_n$, satisfying  $T_{nm}=T_nT_m$. For $x\in X$ and open $U\subseteq X$ put $R(x,U)=\{n\in\mathbb N:T_nx\in U\}$. A set $A\subseteq\mathbb N$ is \emph{dynamically central syndetic in $(\mathbb N,\cdot)$} if $A\supseteq R(x,U)$ for some minimal multiplicative system, some open $U$ and some $x\in U$. A set $A$ is \emph{multiplicatively syndetic} if there is a finite $F\subseteq\mathbb N$ with $\bigcup_{t\in F}\{n: tn\in A\}=\mathbb N$, and \emph{multiplicatively central} if it is central in the semigroup $(\mathbb N,\cdot)$ \cite[Definition 19.20]{HS12}. The family of central sets is denoted by $\cC.$ Similarly, one can define the families $d\cT$, $dc\cT$, $d\cP\cS$, and $dc\cP\cS$, and their respective duals for the multiplicative setting. The following diagram explains their relations:

\[
\xymatrix@C=0.05cm@R=0.1cm{ &&& \cP\cS^*=d \cP\cS^* \ar[ddd] \ar[ldd]\\ &&&&&&\\ &&  \cT \ar[lldd]\ar[rddddd] \\ &&& dc\cP\cS^* \ar[ddr] \ar[ddl]\\ d\cT \ar[ddd] &  \\ && \cC^* \ar[dddddr]\ar[ddll] \ar[ddr] && dc\cS \ar[ddd] \ar[ddl] \\ &&&&&&\\ dc\cT \ar[ddrr]&&& \cC \ar[ddl] \\ &&&  & d\cS \ar[ddl] \\ && dc\cP\cS \ar[dd] \\ &&& \cS \ar[dl] \\ && \cP\cS=d\cP\cS}
\]

We investigate whether dynamically central syndetic sets necessarily contain configurations involving both addition and multiplication. Our main result gives an affirmative answer to this question.

\begin{thma}
    Every dynamically central syndetic subset \(D\) of \(\mathbb N\) contains \(x,y\in\mathbb N\) such that ${\{x,y,x+y,xy\}\subseteq D.}$
    \end{thma}
\begin{proof}
 See Theorem \ref{sum-product}.    
\end{proof}

The theorem establishes that the dynamical structure underlying a dynamically central syndetic set is sufficiently rich to guarantee the simultaneous occurrence of an additive and a multiplicative configuration. 

Further, we show that there is a $dc\mathcal{PS}$ set in $\mathbb{N}$ avoiding the Hindman configurations $\{x,y,xy,x+y\}$ for any $x, y.$

\section{Dynamically Central Syndetic Sets and the Pattern $\{x,y,xy,x+y\}$}

We show that every dynamically central syndetic subset of $(\mathbb N,\cdot)$ contains the configuration $\{x,y,xy,x+y\}$. We use the following result.

\begin{prop}\label{prop3.7}
    Let $P$ be a finite set of integral polynomials without constant coefficients. Any subset of $\mathbb{N}$ that is both multiplicatively syndetic and multiplicatively central contains subsets of the form $\{x,y,xy, x+p(y) : p \in P\}.$
\end{prop}
\begin{proof}
    See \cite[Proposition 3.7]{ABS25}.
\end{proof}

\begin{lemma}\label{mult-dcS}
Every dynamically central syndetic set in $(\mathbb N,\cdot)$ is multiplicatively syndetic and multiplicatively central.
\end{lemma}

\begin{proof}
Let $A$ be a dynamically central syndetic set in $(\mathbb N,\cdot)$. Then by definition, there exists a minimal multiplicative system $(X, (\mathbb N,\cdot))$ with a compact metric space $X$, and there exists an open set $U$ in $X$ and $x\in U$, such that $R(x, U) \subseteq A.$ 

\medskip
We first show the multiplicatively syndeticity of $A$. Since $(X,(\mathbb N,\cdot))$ is minimal, so  for each $y\in X$ the orbit $\{T_ny:n\in\mathbb N\}$ is dense, so $y\in T_t^{-1}U$ for some $t$. Thus $X = \bigcup_{t \in \Lambda} T_t^{-1}U.$ Since $X$ is compact, there exists a finite subset $F\subseteq \Lambda$ such that $X=\bigcup_{t\in F}T_t^{-1}U$. 

\medskip
Given $n\in\mathbb N$, consider the continuous map $T_n: X=\bigcup_{t\in F}T_t^{-1}U \to X=\bigcup_{t\in F}T_t^{-1}U$. There is a $t\in F$ with $T_nx\in T_t^{-1}U$. Since the dynamical system is multiplicative, hence $T_{tn}x\in U$, so $tn\in R(x,U)\subseteq A$, that is, $n \in t^{-1}R(x, U)\subseteq t^{-1}A.$ This proves that $R(x, U)$ (and $A$) is syndetic.

\medskip

Next, we show that $A$ is multiplicatively central. Note that every point of a minimal system is uniformly recurrent \cite[Theorem 1.15]{Fur81}, and $x$ is proximal to itself; since $x\in U$, the dynamical characterization of central sets in countable semigroups \cite[Theorem~19.27]{HS12} shows that $R(x,U)$, hence $A$, is multiplicatively central in $(\mathbb N,\cdot)$.
\end{proof}

\begin{theorem}\label{sum-product}
Let $A\subseteq \mathbb N$ be a dynamically central syndetic set in the semigroup $(\mathbb N,\cdot)$ and let $P$ be a finite set of integral polynomials without constant coefficients. Then there exist $x,y\in \mathbb N$ such that $\{x,y,xy,x+p(y): p \in P\}\subseteq A$.
\end{theorem}

\begin{proof}
By Lemma~\ref{mult-dcS}, $A$ is multiplicatively syndetic and multiplicatively central, and Proposition~\ref{prop3.7} applies.
\end{proof}

\begin{cor}
 Let $A\subseteq \mathbb N$ be a dynamically central syndetic set in the semigroup $(\mathbb N,\cdot)$. Then there exist $x,y\in \mathbb N$ such that $\{x,y,xy,x+y\}\subseteq A$.
\end{cor}

\begin{proof}
    Apply Theorem~\ref{sum-product} with $P=\{p\}$, $p(y)=y$.
\end{proof}

\begin{rem}
    The configuration $\{x,y,xy,x+y\}$ is the one Hindman's conjectured to be partition regular. Theorem~\ref{sum-product} does not bear directly on that conjecture, because $dc\mathcal S$ is not partition regular.
\end{rem}

\begin{prop}\label{rem:notPR}
The family $dc\mathcal S$ is not partition regular in $(\mathbb N,+)$ or in $(\mathbb N,\cdot)$.
\end{prop}

\begin{proof}
We first show that the family $dc\cS$ is not partition regular in $(\mathbb N, +)$. To prove this, let us consider a $2$-colouring on $(\mathbb{N}, +)$ given  by 
\[
c: (\mathbb{N}, +) \to \{0,1\}\qquad n \mapsto \lceil \log_2(n)\rceil \pmod 2.
\]
\noindent For this colouring, both the pieces $c^{-1}(0)=\displaystyle\bigsqcup_{m\ge 0}[2^{2m, 2^{2m+1}})$ and  $c^{-1}(1)=\displaystyle \bigsqcup_{m\ge 0}[2^{2m+1}, 2^{2m+2})$ are not syndetic.

Next for the remaining part, we  consider the semigroup homomorphism $\Omega: (\mathbb{N}, .) \to (\mathbb{N}, +)$, where $\Omega(n)$ the number of prime factors of $n$ counted with multiplicity. Then the composition
\[
\bar{c}:=c\circ \Omega : (\mathbb{N}, .)  \to \{0,1\} \qquad n \mapsto  \lceil \log_2(\Omega(n))\rceil \pmod 2
\]
gives a coloring of $(\mathbb{N}, .).$ Using a similar argument, one may prove that none of the pieces $\bar{c}^{-1}(0)$ and $\bar{c}^{-1}(1)$ are syndetic. 
\end{proof}

The hypothesis of Theorem~\ref{sum-product} cannot be weakened to the corresponding piecewise syndetic notion, even if one asks for largeness in both semigroup structures at once.

\begin{defn}
Let $B_1,B_2,\dots$ be finite subsets of $\{2,3,\dots\}$ with $\min B_{s+1}>(\max B_s)^2$. For odd $k$ let $B_k=[N_k,N_k+k]$ with $N_k>k$. For even $k$ let $B_k=q_k\cdot\{1,\dots,k\}$ with $q_k>k$ prime. Put $A=\bigcup_kB_k$.
\end{defn}

\begin{prop}
    The set $A$ is thick in both  $(\mathbb N,+)$ and  $(\mathbb N,\cdot)$.
\end{prop}
\begin{proof}
    The first holds because of the odd blocks. For the second, given a finite $F$, take an even $k\ge\max F$; then $F\cdot q_k\subseteq A$. A thick set meets every syndetic set, so in either semigroup $A$ meets every dynamically central syndetic set; that is, $A$ is a set of pointwise recurrence. Hence $A=\mathbb N\cap A$ lies in $dc\mathcal{PS}$ in both $(\mathbb N,+)$ and $(\mathbb N,\cdot)$. It is also central in both, since thick sets are central \cite[Theorem~4.48]{HS12}.
\end{proof}

Now we are in a position to show that the above constructed set $A$ does not contain any Hindman configuration.

\begin{cor}\label{thickavoid}
 The set $A$ contains no configuration $\{x,y,xy,x+y\}$.  
\end{cor}
\begin{proof}
    Suppose $2\le x\le y$ with $y\in B_s$ and $x+y,xy\in A$. Then $y<x+y$ and $y\le xy$, and both are at most $(\max B_s)^2<\min B_{s+1}$. Hence $x+y,xy\in B_s$.
    
\noindent {\bf Case I.}~If $B_s=[N,N+k]$, then $xy\ge2N>N+k$. Thus $xy$ does not belong to $B_s$. This is a contradiction.

\noindent {\bf Case II.}~If $B_s=q\cdot\{1,\dots,k\}$, then $q$ divides $y$ and $x+y$, hence $x$. Every element of an earlier block is smaller than $q$, so $x\in B_s$. Writing $x=qa$ and $y=qb$, we get $xy=q\cdot(qab)$ with $qab\ge q>k$, so $xy\notin B_s$. This is also a contradiction.

Hence we derive that the set $A$ does not contain the pattern $\{x,y,xy,x+y\}$ for any choice of $x, y.$
\end{proof}

\vspace{1cm}

\end{document}